\documentclass[twoside,a4paper,12pt,centertags,reqno]{amsart} % 13pt displays better than 12pt
\usepackage{amsmath,amssymb,verbatim,vmargin}
\usepackage{color}
\usepackage{tikz}
\usepackage{hyperref}
\hypersetup{colorlinks,bookmarks=true,linktocpage=true,citecolor=blue,linkcolor=magenta}

\usepackage{mathptmx}   % based on the Times Roman font

\allowdisplaybreaks % align break
\usepackage{color}
\usepackage[all]{xy} % diagrams

\usepackage{mathtools}
\usepackage{MnSymbol} %wedge righthalfcup

\DeclareFontFamily{U}{mathx}{\hyphenchar\font45}
\DeclareFontShape{U}{mathx}{m}{n}{
      <5> <6> <7> <8> <9> <10>
      <10.95> <12> <14.4> <17.28> <20.74> <24.88>
      mathx10
      }{}
\DeclareSymbolFont{mathx}{U}{mathx}{m}{n}
\DeclareFontSubstitution{U}{mathx}{m}{n}
\DeclareMathAccent{\widecheck}{0}{mathx}{"71}
\DeclareMathAccent{\wideparen}{0}{mathx}{"75}

\theoremstyle{plain}
\newtheorem{thm}{Theorem}[section]

\theoremstyle{definition}

\newtheorem{rem}[thm]{Remark}

\usepackage{enumerate}

\newcommand{\pd}{\partial}

\newcommand{\bR}{{\mathbb R}}

\newcommand{\diam}{\text{{\rm diam}}\,}
\def\barint_#1{\mathchoice
            {\mathop{\vrule width 6pt
height 3 pt depth -2.5pt
                    \kern -9.5pt
\intop \kern -4pt}\nolimits_{#1}}%
            {\mathop{\vrule width 5pt height
3 pt depth -2.6pt
                    \kern -6.5pt
\intop \kern -4pt}\nolimits_{#1}}%
            {\mathop{\vrule width 5pt height
3 pt depth -2.6pt
                    \kern -6pt
\intop \kern -4pt}\nolimits_{#1}}%
            {\mathop{\vrule width 5pt height
3 pt depth -2.6pt
          \kern -6pt \intop \kern -4pt}\nolimits_{#1}}}
          
           \def\bariint_#1{\mathchoice
            {\mathop{\vrule width 15pt
height 3 pt depth -2.5pt
                    \kern -15.8pt
\intop \kern -8pt\intop \kern -4pt}\nolimits_{#1}}%
            {\mathop{\vrule width 9pt height
3 pt depth -2.6pt
                    \kern -10.5pt
\intop \kern -8pt\intop \kern -4pt}\nolimits_{#1}}%
            {\mathop{\vrule width 9pt height
3 pt depth -2.6pt
                    \kern -10pt
\intop \kern -8pt\intop \kern -4pt}\nolimits_{#1}}%
            {\mathop{\vrule width 9pt height
3 pt depth -2.6pt
          \kern -8pt \intop \kern -10pt\intop \kern -4pt}
      \nolimits_{  #1}}}

\def\barintlim_#1{\mathchoice
            {\mathop{\vrule width 6pt
height 3 pt depth -2.5pt
                    \kern -8.8pt
\intop \kern -4pt}\limits_{#1}}%
            {\mathop{\vrule width 5pt height
3 pt depth -2.6pt
                    \kern -6.5pt
\intop \kern -4pt}\limits_{#1}}%
            {\mathop{\vrule width 5pt height
3 pt depth -2.6pt
                    \kern -6pt
\intop \kern -4pt}\limits_{#1}}%
            {\mathop{\vrule width 5pt height
3 pt depth -2.6pt
          \kern -6pt \intop \kern -4pt}\limits_{#1}}}
          
           \def\bariintlim_#1{\mathchoice
            {\mathop{\vrule width 15pt
height 3 pt depth -2.5pt
                    \kern -15.8pt
\intop \kern -8pt\intop \kern -4pt}\limits_{#1}}%
            {\mathop{\vrule width 9pt height
3 pt depth -2.6pt
                    \kern -10.5pt
\intop \kern -8pt\intop \kern -4pt}\limits_{#1}}%
            {\mathop{\vrule width 9pt height
3 pt depth -2.6pt
                    \kern -10pt
\intop \kern -8pt\intop \kern -4pt}\limits_{#1}}%
            {\mathop{\vrule width 9pt height
3 pt depth -2.6pt
          \kern -8pt \intop \kern -10pt\intop \kern -4pt}
      \limits_{  #1}}}
          
\renewcommand{\iint}{\int \kern -3pt\int}       

\numberwithin{equation}{section}
\makeatletter
\@namedef{subjclassname@2020}{\textup{2020} Mathematics Subject Classification}
\makeatother

\newcommand{\nocontentsline}[3]{}
\let\origcontentsline\addcontentsline
\newcommand\stoptoc{\let\addcontentsline\nocontentsline}
\newcommand\resumetoc{\let\addcontentsline\origcontentsline}

\title{Refined Qualitative Properties of Euclidean Heat Content}

\author{Yi C. Huang} 
\address{School of Mathematical Sciences, Nanjing Normal University, Nanjing 210023, People's Republic of China}
\email{Yi.Huang.Analysis@gmail.com}
\urladdr{https://orcid.org/0000-0002-1297-7674}

\date{\today} 

\thanks{Research of the Author is partially supported by the National NSF grant of China (no. 11801274), the JSPS Invitational Fellowship for Research in Japan (no. S24040),
and the Open Projects from Yunnan Normal University (no. YNNUMA2403) and Soochow University (no. SDGC2418). 
The Author thanks in particular Professors Michiel van den Berg, Xueping Huang and Tohru Ozawa for helpful communications.}

\begin{document}

\begin{abstract}
Invoking a rather simple semigroup approach, we refine the recent qualitative properties of Euclidean heat content obtained by van den Berg and Gittins.
\end{abstract}

\maketitle

%\tableofcontents

%\section{Introduction}

In this paper we refine quantitatively the qualitative properties (in particular, the convexity and monotonicity obtained by van den Berg and Gittins in \cite{berg2024qualitative}) of heat content
%\begin{equation}  
$$H_\Omega(t)=\int_\Omega\int_\Omega p_t(x,y)dxdy,\quad t\geq0,$$
%\end{equation}
where $\Omega$ is a non-empty open set in $\bR^m$ and 
%\begin{equation}  
$$p_t(x,y)=\frac{1}{(4\pi t)^{\frac{m}{2}}}e^{-\frac{|x-y|^2}{4t}}$$
%\end{equation}
is the standard heat kernel in $\bR^m$. Using the heat semigroup 
%\begin{equation}  
$$e^{t\Delta}f(x)=\int_{\bR^m}p_t(x,y)f(y)dy$$
%\end{equation}
and the $L^2(\bR^m)$ pairing for real valued functions
%\begin{equation}  
$$(f,g)=\int_{\bR^m}f(x)g(x)dx,$$
%\end{equation}
we can rewrite the heat content into the following convenient form
\begin{equation} \label{e:form}
H_\Omega(2t)=(e^{2t\Delta}\chi_\Omega,\chi_\Omega)=(e^{t\Delta}\chi_\Omega,e^{t\Delta}\chi_\Omega).
\end{equation}
Via this formulation it is immediate that the heat content $H_\Omega(t)$ is decreasing and convex.

\bigskip

We can also express the second order derivative of heat content in a transparent way
\begin{equation} \label{e:convexkey}
\frac{d^2}{dt^2}H_\Omega(2t)=4(\Delta e^{t\Delta}\chi_\Omega, \Delta e^{t\Delta}\chi_\Omega).
\end{equation}
Notice that
%\begin{equation}
$$\Delta (e^{t\Delta}\chi_\Omega)(x)=\int_\Omega\pd_tp_t(x,y)dy,$$
%\end{equation}
whereas the time derivative of heat kernel satisfies
\begin{equation} \label{e:kernel}
\pd_tp_t(x,y)=\frac1t \left[-\frac{m}{2}+\frac{|x-y|^2}{4t}\right]p_t(x,y).
\end{equation}
Suppose $\diam(\Omega)<\infty$. Therefore, for $t\geq(\diam(\Omega))^2$, as $|x-y|\leq \diam(\Omega)$,
\begin{equation} \label{e:key}
\pd_tp_t(x,y)\leq -\frac{2m-1}{4}\frac{p_t(x,y)}{t}.
\end{equation}
Hence, integrating \eqref{e:key} with respect to $(x,y)\in\Omega\times\Omega$ gives
\begin{equation} \label{e:firstorderineq}
\frac{d}{dt}H_\Omega(t)\leq -\frac{2m-1}{4}\frac{H_\Omega(t)}{t},\quad t\geq{(\diam(\Omega))^2}.
\end{equation}

\begin{rem}
\eqref{e:firstorderineq} improves the monotonicity result of van den Berg and Gittins \cite{berg2024qualitative}
%\begin{equation}
$$\frac{d}{dt}H_\Omega(t)\leq -\frac{4m^2+4m-7}{8(m+2)e^{\frac14}}\frac{H_\Omega(t)}{t},\quad t\geq {(\diam(\Omega))^2}.$$
%\end{equation}
\end{rem}

\bigskip

Meanwhile, inserting \eqref{e:key} into \eqref{e:convexkey}, we get
%\begin{equation}
$$\frac{d^2}{dt^2}H_\Omega(2t)\geq 4\left(\frac{2m-1}{4}\right)^2\frac{H_\Omega(2t)}{t^2},\quad t\geq(\diam(\Omega))^2,$$
%\end{equation}
Equivalently, this can be written as
\begin{equation} \label{e:main1}
\frac{d^2}{dt^2}H_\Omega(t)\geq \left(\frac{2m-1}{2}\right)^2\frac{H_\Omega(t)}{t^2},\quad t\geq2{(\diam(\Omega))^2}.
\end{equation}

\begin{rem}
\eqref{e:main1} improves, with essentially the same large-time condition, the convexity result of van den Berg and Gittins \cite{berg2024qualitative} 
(they considered the $\pd^2_{tt}p_t(x,y)$ version of \eqref{e:kernel})
%\begin{equation}
$$\frac{d^2}{dt^2}H_\Omega(t)\geq \frac{4m^2+4m-7}{16}\frac{H_\Omega(t)}{t^2},\quad t\geq{(\diam(\Omega))^2}.$$
%\end{equation}
\end{rem}

\bigskip

Furthermore, using the heat kernel lower bound
%\begin{equation}
$$p_t(x,y)\geq\frac{e^{-\frac18}}{(4\pi t)^{\frac{m}{2}}},\quad t\geq2{(\diam(\Omega))^2},$$
%\end{equation}
we get
\begin{equation} \label{e:lowerHeatContent}
H_\Omega(t)\geq e^{-\frac18}\frac{|\Omega|^2}{(4\pi t)^{\frac{m}{2}}},
\end{equation}
which in conjunction with \eqref{e:main1} gives
\begin{equation} \label{e:main2}
\frac{d^2}{dt^2}H_\Omega(t)\geq \left(\frac{2m-1}{2}\right)^2\frac{e^{-\frac18}}{(4\pi)^{\frac{m}{2}}}\frac{|\Omega|^2}{t^{\frac{m}{2}+2}},\quad t\geq2{(\diam(\Omega))^2}.
\end{equation}
Integrating this inequality over the interval $[t,+\infty)$, with $t\geq 2{(\diam(\Omega))^2}$, gives
\begin{equation} \label{e:main3}
\frac{d}{dt}H_\Omega(t)\leq -\frac{1}{\frac{m}{2}+1}\left(\frac{2m-1}{2}\right)^2\frac{e^{-\frac18}}{(4\pi)^{\frac{m}{2}}}\frac{|\Omega|^2}{t^{\frac{m}{2}+1}},\quad t\geq2{(\diam(\Omega))^2}.
\end{equation}
This (except for $m=1$) is sharper than the resulted estimate by merely combining \eqref{e:firstorderineq} and \eqref{e:lowerHeatContent}.
Since $H_\Omega(t)$ is convex, its first derivative is increasing and continuous, hence
\begin{equation} \label{e:main3small}
\frac{d}{dt}H_\Omega(t)\leq -\frac{1}{\frac{m}{2}+1}\left(\frac{2m-1}{2}\right)^2\frac{e^{-\frac18}}{(4\pi)^{\frac{m}{2}}}\frac{|\Omega|^2}{(2{(\diam(\Omega))^2})^{\frac{m}{2}+1}},\quad t\leq2{(\diam(\Omega))^2}.
\end{equation}
By similar reasoning (i.e., $\frac{d^3}{dt^3}H_\Omega(t)\leq0$, verified again via \eqref{e:form}), we deduce from \eqref{e:main2}
\begin{equation} \label{e:main2small}
\frac{d^2}{dt^2}H_\Omega(t)\geq \left(\frac{2m-1}{2}\right)^2\frac{e^{-\frac18}}{(4\pi)^{\frac{m}{2}}}\frac{|\Omega|^2}{(2{(\diam(\Omega))^2})^{\frac{m}{2}+2}},\quad t\leq2{(\diam(\Omega))^2}.
\end{equation}

\begin{rem}
In the small-time case for $\bR^m$, there is no strict-convexity result in \cite{berg2024qualitative}.
\end{rem}

\bigskip

In summary, we have obtained the following quantitatively refined properties (listed in the order of \cite[Theorem 3]{berg2024qualitative}) of the heat content: 
\eqref{e:main1}, \eqref{e:main2}\&\eqref{e:main2small}, \eqref{e:main3}\&\eqref{e:main3small}, \eqref{e:firstorderineq}.
The novelty lies in using the semigroup formulation \eqref{e:form} of the heat content.

\bigskip

\stoptoc

\section*{\textbf{Compliance with ethical standards}}

\bigskip

\textbf{Conflict of interest} The author has no known competing financial interests
or personal relationships that could have appeared to influence this reported work.

\bigskip

\textbf{Availability of data and material} Not applicable.

\bigskip

\resumetoc

\bibliographystyle{alpha}

\bibliography{Hua-RefinedHeatContent} 
 
\end{document}